# ON THE BURNSIDE PROBLEM ON PERIODIC GROUPS

Sergei V. Ivanov

ABSTRACT. It is proved that the free $m$-generated Burnside groups $\mathbb{B}(m, n)$ of exponent $n$ are infinite provided that $m > 1$, $n \geq 2^{48}$.

In 1902 William Burnside posed the following problem [2]. *Does a group $G$ have to be finite provided that $G$ has a finite set of generators and its elements satisfy the identity $x^n = 1$?* In other words, must a finitely generated group $G$ of exponent $n$ be finite?

In the same paper, Burnside proved that the problem was solved in the affirmative for groups of exponents 2, 3 and for 2-generated groups of exponent 4 as well.

In 1940 Sanov [12] obtained a positive solution to the Burnside problem for the case of exponent 4.

The next significant step was made by Marshall Hall [4] in 1957 when he solved the problem in the affirmative for the exponent of 6.

In 1964 Golod [3] found the first example of an infinite periodic group with a finite number of generators. Although that example did not satisfy the identity $x^n = 1$, i.e., the group was of unbounded exponent, it gave the first positive evidence that the Burnside problem might not be solved affirmatively for all exponents (and it might possibly fail for very large exponents).

In 1968 Novikov and Adian achieved a real breakthrough in a series of fundamental papers [9] in which some ideas put forward by Novikov [8] in 1959 were developed to prove that there are infinite periodic groups of odd exponents $n \geq 4381$ with $m > 1$ generators. Later, Adian [1] improved the estimate up to $n \geq 665$ ($n$ is odd again). Notice in the papers [9] that, in fact, the free Burnside groups $\mathbb{B}(m, n) = \mathbb{F}_m / \mathbb{F}_m^n$, where $\mathbb{F}_m$ is a free group of rank $m > 1$ and $\mathbb{F}_m^n$ is the normal subgroup of $\mathbb{F}_m$ generated by all $n$th powers (with odd $n \geq 4381$) of elements of $\mathbb{F}_m$, were constructed and studied. Using a very complicated inductive construction, Novikov and Adian presented the group $\mathbb{B}(m, n)$ by defining relations of the

1991 *Mathematics Subject Classification.* Primary 20F05, 20F06, 20F32, 20F50.

Received by the editors January 7, 1992

Lectures on results of this note were given at the University of Utah (October 30, 1991, January 9, 16, 23, 30, 1992), the University of Wisconsin-Parkside (February 6, 1992), the City University of New York (February 7, 1992), the University of Nebraska-Lincoln (February 13, 14, 1992), the Kent State University (March 2, 1992), the University of Florida-Gainesville (March 23, 1992). It is the pleasure of the author to thank these universities for sponsoring his visits as well as to express his gratitude to Professors G. Baumslag, B. Chandler, P. Enflo, S. Gagola, S. Gersten, J. Keesling, A. Lichtman, S. Margolis, J. Meakin, G. Robinson, and J. Thompson for their interest in this work form $A^n = 1$, where $A$'s are some specially chosen elements of $\mathbb{F}_m$, and studied their consequences. They not only obtained the result that the group $\mathbb{B}(m, n)$ is infinite but also other important information about $\mathbb{B}(m, n)$. For example, it was proved that the word and conjugacy problems are solvable in $\mathbb{B}(m, n)$ and that any finite or abelian subgroup of $\mathbb{B}(m, n)$ is cyclic (under the restrictions on $m$ and $n$ above; for these and other results see [1]).

At the same time, it should be pointed out that [9] is very long and of very complicated logical structure.

In 1982 Ol′shanskiĭ [10] succeeded in finding a considerably shorter proof of the theorem of Novikov and Adian, although the estimate $n > 10^{10}$ (where $n$ again is odd) of [10] is much worse than $n \geq 665$ of Adian's [1]. On the other hand, it is worth noting that the approach of Ol′shanskiĭ's to treat the free Burnside groups $\mathbb{B}(m, n)$ is based on a powerful geometric method of graded diagrams (see [11, 6] for numerous applications of the method in combinatorial group theory).

Thus, it is known that the Burnside problem is settled in the affirmative for exponents $n = 2, 3, 4, 6$ and in the negative for the exponents that have an odd divisor not less than 665 (the latter is an easy corollary of the theorem of Novikov and Adian). In particular, the Burnside problem still remains open for exponents of the form $n = 2^k$. Besides, there is no approach to study the free Burnside groups $\mathbb{B}(m, n)$ of even exponent $n$, even if $n$ has a rather great odd divisor, and the only known characteristic of these groups is their infiniteness (unlike the case of odd exponents $n \geq 665$).

Now let us mention an unpublished work [5] of the author's where the free Burnside groups $\mathbb{B}(m, 2n)$ and $\mathbb{B}(m, 4n)$ with odd $n \gg 1$ were constructed by means of defining relations in order to prove solvability of the word and conjugacy problems for these groups and to obtain a description of their finite subgroups. All efforts to extend the techniques of [5] to study the groups $\mathbb{B}(m, 8n)$ with odd $n \gg 1$, however, were unsuccessful.

In the meantime, quite new techniques have been developed in order to construct and study the free Burnside groups $\mathbb{B}(m, n)$ with any $n \gg 1$ regardless of the oddness of $n$. The key point of the techniques is in obtaining a complete description of finite subgroups of the free Burnside groups $\mathbb{B}(m, n)$. Therefore, in Theorem A, which gives the negative solution to the problem of Burnside's for all rather great exponents, we include this description.

**Theorem A.** *Let $\mathbb{B}(m, n)$ be the free Burnside group of rank $m$ and exponent $n$, $m > 1$ and $n \geq 2^{48}$. Then*

(a) *The group $\mathbb{B}(m, n)$ is infinite.*

(b) *The word and conjugacy problems are solvable in $\mathbb{B}(m, n)$.*

(c) *Suppose $n = 2^k n_0$, where $n_0$ is odd. If $k = 0$ (i.e. $n$ is odd) then any finite subgroup of $\mathbb{B}(m, n)$ is cyclic. If $k > 0$ (i.e. $n$ is even) then any finite subgroup of $\mathbb{B}(m, n)$ is isomorphic to a subgroup of a direct product of two groups, one of which is a dihedral group of order $2n$ and the other is a direct product of several copies of a dihedral group of order $2^{k+1}$. In particular, if $n = 2^k$ then any finite subgroup of $\mathbb{B}(m, n)$ is just a subgroup of a direct product of several copies of a dihedral group of order $2n$.*

(d) *The center of the group $\mathbb{B}(m, n)$ is trivial.*

Now let us give an inductive construction of the group $\mathbb{B}(m, n)$ of any exponent $n \gg 1$ by means by defining relations. Notice that this construction repeats (it is a surprise in itself!) a construction invented by Ol'shanskiĭ [10] for the case where $n$ is odd.

On the set of all nonempty reduced words over an alphabet $\mathbb{A} = \{a_1^{\pm 1}, \ldots, a_m^{\pm 1}\}$ (we assume $\mathbb{F}_m$ to be the free group over the alphabet $\mathbb{A}$), we introduce a total order $a_1 \prec a_2 \prec \cdots$ such that $|X| \leq |Y|$ implies $X \preceq Y$, where $|X|$ denotes the length of the word $X$.

Now, for each $i \geq 1$, we define a word $A_i$ called the period of rank $i$ to be the smallest (in terms of the order "$\prec$" introduced above) of those words over $\mathbb{A}$ whose orders in the group $\mathbb{B}(i-1)$, given by the presentation

$$(*) \qquad \mathbb{B}(i-1) = \langle a_1^{\pm 1}, \ldots, a_m^{\pm 1} \mid\mid A_1^n = 1, \ldots, A_{i-1}^n = 1 \rangle,$$

are infinite.

Notice that it is not clear a priori whether $A_i$ exists for each $i$ or not. Notice also that infiniteness of the free Burnside groups $\mathbb{B}(m,n)$ (under the restrictions on $m$, $n$ above) follows from the next Theorem B, since a finite group cannot be presented by infinitely many independent defining relations over a finite alphabet.

**Theorem B.** *Suppose $m > 1$ and $n \geq 2^{48}$. Then the period $A_i$ of rank $i$ does exist for each $i \geq 1$, i.e., the system $\{A_i^n = 1\}_{i=1}^{\infty}$ is infinite. Next, the system $\{A_i^n = 1\}_{i=1}^{\infty}$ can be taken as an independent set of defining relations of the free Burnside group $\mathbb{B}(m,n)$ and order of the period $A_i$ of any rank $i \geq 1$ is equal in $\mathbb{B}(m,n)$ to $n$ exactly.*

The following theorem contains some basic technical results about finite subgroups of the groups $\mathbb{B}(i-1)$ and $\mathbb{B}(m,n)$. Notice that one can derive the algebraic description of finite subgroups of $\mathbb{B}(m,n)$ given in Theorem A (proceeding by induction on the maximum of heights of words of a finite subgroup of $\mathbb{B}(m,n)$) from (a)–(e) of the following.

**Theorem C.** *Let $\mathbb{B}(m,n)$ be the free Burnside group of rank $m > 1$ and exponent $n \geq 2^{48}$, and suppose that $\mathcal{F}(A_i)$ is a maximal finite subgroup of the group $\mathbb{B}(i-1)$ given by $(*)$ with respect to the property that $\mathcal{F}(A_i)$ is normalized by the period $A_i$ of rank $i$. Next, denote by $J_i$ a word such that the inclusions $J_i^2, (J_i A_i)^2 \in \mathcal{F}(A_i)$, hold in $\mathbb{B}(i-1)$ (if there exists no such word we simply put $J_i = 1$). Then the following claims hold*:

(a) *Any word $W$ having finite order in $\mathbb{B}(i-1)$ is conjugate in $\mathbb{B}(i-1)$ to a word of the form $A_j^k T$ for some integer $k$, $j < i$ and $T \in \mathcal{F}(A_j)$. Besides, conjugacy in $\mathbb{B}(i-1)$ of the words $A_{j_1}^{k_1} T_1$ and $A_{j_2}^{k_2} T_2$, where $T_1 \in \mathcal{F}(A_{j_1})$ and $T_2 \in \mathcal{F}(A_{j_2})$, $j_1, j_2 < i$, that are not equal in $\mathbb{B}(i-1)$ to the identity yields $j_1 = j_2$ and $k_1 = \pm k_2 \pmod{n}$. (Therefore, given a nontrivial word $W$, such a number $j$ is defined uniquely in $\mathbb{B}(m,n)$ as well as in $\mathbb{B}(i-1)$ and called the height of the word $W$.)*

(b) *$\mathcal{F}(A_i)$ is defined uniquely, embeds into $\mathbb{B}(m,n)$, consists of words whose heights are less than $i$, and is a 2-group.*

(c) *Any finite subgroup of $\mathbb{B}(m,n)$ consisting of words of heights $\leq i$ and containing a word of height $i$ exactly is conjugate to a subgroup of the group generated by $A_i$, $J_i$, and all words from $\mathcal{F}(A_i)$.*

(d) *The subgroup $\mathcal{F}(A_i)$ of $\mathbb{B}(i-1)$ is normalized by $J_i$.*

(e) *The words $J_i$ and $A_i$ act on the subgroup $\mathcal{F}(A_i)$ of $\mathbb{B}(i-1)$ by conjugations in the same way as some words $V_1$ and $V_2$ act respectively, where $V_1$ and $V_2$ are* such that the subgroup of $\mathbb{B}(i-1)$ generated by $V_1$, $V_2$ and by all words from $\mathcal{F}(A_i)$ is finite and the equation $J_i^2 = V_1^2$ (as well as $(J_i A_i)^2 = (V_1 V_2)^2$ provided $J_i \neq 1$) holds in $\mathbb{B}(i-1)$.

Let us conclude with some remarks about proofs of Theorems A, B, and C.

First, in the case of odd $n$ (this special case emerges as the simplest one where the finite subgroups $\mathcal{F}(A_i)$ are trivial for all $i$), proofs of Theorems A, B, and C virtually repeat the proofs of Ol′shanskiĭ's [10]. In particular, we use a geometric interpretation of deducibility of relations in a group from its defining relations (this interpretation is based on the notion of van Kampen diagrams, see [7]).

On the other hand, the case where $n$ is even requires much more delicate investigations of various properties of finite subgroups of the free Burnside group $\mathbb{B}(m,n)$. As a matter of fact, we point out that these properties of finite subgroups of groups $\mathbb{B}(m,n)$ along with subgroups $\mathcal{F}(A_i)$ degenerate in the case of odd $n$, and so one can say that the works [1, 9, 10] primarily deal with most general characteristics of the groups $\mathbb{B}(m,n)$.

Finally, we mention that our estimate $n \geq 2^{48}$ is rather rough and can be strongly improved at cost of complication of proofs.

*Added in proof*. It has been known to the author that I. Lysionak, *The infinity of Burnside groups of exponents $2^\kappa$ for $\kappa \geq 13$* (preprint), announces an independent solution of the Burnside problem for exponents of the form $2^\kappa \geq 2^{13}$ based on the Novikov-Adian method.

## Acknowledgments

The author is grateful to Professors Steve Gersten and Alexander Ol′shanskiĭ for helpful discussions and their encouragement.## References

1. S. I. Adian, *The Burnside problems and identities in groups*, Moscow, Nauka, 1975.
2. W. Burnside, *On unsettled question in the theory of discontinuous groups*, Quart. J. Pure Appl. Math. **33** (1902), 230–238.
3. E. S. Golod, *On nil-algebras and finitely residual groups*, Izv. Akad. Nauk SSSR. Ser. Mat. **28** (1964), 273–276.
4. M. Hall, *Solution of the Burnside problem for exponent 6*, Proc. Nat. Acad. Sci. U.S.A. **43** (1957), 751–753.
5. S. V. Ivanov, *Free Burnside groups of some even exponents*, 1987 (unpublished).
6. S. V. Ivanov and A. Yu. Ol′shanskiĭ, *Some applications of graded diagrams in combinatorial group theory*, London Math. Soc. Lecture Note Ser., vol. 160 (1991), Cambridge Univ. Press, Cambridge and New York, 1991, 258–308.
7. R. C. Lyndon and P. C. Schupp, *Combinatorial group theory*, Springer-Verlag, Heidelberg, 1977.
8. P. S. Novikov, *On periodic groups*, Dokl. Akad. Nauk SSSR Ser. Mat. **27** (1959), 749–752.
9. P. S. Novikov and S. I. Adian, *On infinite periodic groups* I, II, III, Izv. Akad. Nauk SSSR. Ser. Mat. **32** (1968), 212–244; 251–524; 709–731.
10. A. Yu. Ol′shanskiĭ, *On the Novikov-Adian theorem*, Mat. Sb. **118** (1982), 203–235.
11. \_\_\_\_\_\_, *Geometry of defining relations in groups*, Moscow, Nauka, 1989.
12. I. N. Sanov, *Solution of the Burnside problem for exponent 4*, Uchen. Zap. Leningrad State Univ. Ser. Mat. **10** (1940), 166–170.Higher Algebra, Department of Mathematics, Moscow State University, Moscow 119899 Russia

*Current address*: Department of Mathematics, University of Utah, Salt Lake City, Utah 84112

*E-mail address*: `ivanov@math.utah.edu`